\documentclass[11pt]{amsart}
\usepackage{latexsym,amssymb,amsmath,amscd,amsthm}
\numberwithin{equation}{section}
\theoremstyle{remark}

\newtheorem{corollary}{{\bf COROLLARY}}[section]

\newtheorem{proposition}{{\bf PROPOSITION}}[section]
\newcommand{\bq}{\begin{equation}}
\newcommand{\bea}{\begin{array}}
\newcommand{\eea}{\end{array}}

\newcommand{\ga}{\alpha}

\newcommand{\gD}{\Delta}
\newcommand{\gl}{\lambda}

\newcommand{\gb}{\beta}

\newcommand{\mf}{\mathfrak}

\newcommand{\wg}{\wedge}

\newcommand{\ul}[1]{\underline{#1}}

\newcommand{\ola}{\overleftarrow}
\newcommand{\go}{\omega}
\newcommand{\gO}{\Omega}
\newcommand{\gG}{\Gamma}

\newcommand{\gz}{\zeta}
\newcommand{\gag}{\gamma}
\newcommand{\gd}{\delta}
\newcommand{\pp}{\partial}

\newcommand{\ora}{\overrightarrow}

\newcommand{\tl}{\tilde}
\newcommand{\na}{\nabla}
\newcommand{\gk}{\kappa}

\newcommand{\ncint}{-\!\!\!\!\!\int\!\!\!\!-}

\title{DISCRETIZATION AND MOYAL BRACKETS}
\author{Robert Carroll\\University of Illinois, Urbana, IL 61801}

\date{January, 2001 - email: rcarroll@math.uiuc.edu}

\begin{document}

\bibliographystyle{plain}

\begin{abstract} {We give a q-analysis version of a discretizaton procedure of
Kemmoku and Saito leading to an apparently new q-Moyal type bracket}
\end{abstract}

\maketitle




\section{INTRODUCTION}
\renewcommand{\theequation}{1.\arabic{equation}}
\setcounter{equation}{0}

We are pursuing further some of the directions 
spelled out in \cite{ch} relating Moyal-Weyl-Wigner
theory, Hirota formulas, integrable systems, and discretization, with additional
connections involving quantum groups (cf. \cite{ca,czd}).  In this note we indicate
an apparently new q-Moyal type bracket formula arising in this context.
In particular we follow here
frameworks from
\cite{ch,kas,kaw}
for deformation quantization and integrable systems and refer to \cite{ca,czd} and
references cited there for q-analysis and quantum groups.
One objective will be to examine various formulas arising in the deformation of
integrable systems and see if there are quantum group versions.
Further we are looking for q-analysis versions of deformation quantization
formulas in order to compare q-calculus and quantum group theory with deformation
quantization.  Thus for background one recalls for wave functions
$\psi$ there are Wigner functions {\bf (WF)} given via
\bq
f(x,p)=\frac{1}{2\pi}\int dy\psi^*\left(x-\frac{\hbar}{2}y\right)exp(-iyp)
\psi\left(x+\frac{\hbar}{2}y\right)
\label{ZZ1}
\end{equation}
Then defining $f*g$ via 
\bq
f*g=f\, exp\left[\frac{i\hbar}{2}(\overleftarrow{\partial}_x
\overrightarrow{\partial}_p-\overleftarrow{\partial}_p\overrightarrow{\partial}_x
\right]\,g;
\label{ZZ2}
\end{equation}
$$f(x,p) * g(x,p)=f\left(x+\frac{i\hbar}{2}\overrightarrow{\partial}_p,
p-\frac{i\hbar}{2}\overrightarrow{\partial}_x\right)g(x,p)$$
time dependence of WF's is given by ($H\sim$ Hamiltonian)
\bq
\partial_tf(x,p,t)=\frac{1}{i\hbar}(H * f(x,p,t)-f(x,p,t) * H)=
\{H,f\}_M
\label{ZZ3}
\end{equation}
where $\{f,g\}_M\sim$ Moyal bracket.  As $\hbar\to 0$ this reduces to $\pp_tf-
\{H,f\}=0$ (standard Poisson bracket).  One can generalize and write out 
\eqref{ZZ2} in various ways.  For example replacing $i\hbar/2$ by $\gk$ one obtains as
in \cite{gal}
\bq
f*g=\sum_0^{\infty}\frac{\kappa^s}{s!}\sum_{j=0}^s(-1)^j{s \choose j}
(\partial_x^j\partial_p^{s-j}f)(\partial_x^{s-j}\partial_p^jg)
\label{X17}
\end{equation}
leading to ($\{f,g\}_{\gk}=(f*g-g*f)/2\gk$)
\bq
\{f,g\}_{\kappa}=\sum_0^{\infty}\frac{\kappa^{2s}}{(2s+1)!}\sum_{j=0}^{2s+1}
(-1)^j{2s+1 \choose j}(\partial_x^j\partial_p^{2s+1-j}f)
(\partial_x^{2s+1-j}\partial_p^jg)
\label{X18}
\end{equation}
(cf. also \cite{sh}) which will also be utilized in the form
\bq\label{1}
f*g=fe^{\gk(\ola{\pp}_x\ora{\pp}_p-\ola{\pp}_p\ora{\pp}_x)}g=e^{[\gk(\pp_{x_1}\pp_{p_2}
-\pp_{x_2}\pp_{p_1})]}f(x_1,p_1)g(x_2,p_2)|_{(x,p)}=
\end{equation}
$$=\sum_0^{\infty}\frac{(-1)^r\gk^{r+s}}{r!s!}\frac{\pp^{r+s}f}{\pp x^r\pp p^s}\frac
{\pp^{r+s}g}{\pp p^r\pp x^s}=\sum_0^{\infty}\frac{\gk^n(-1)^{n-s}}{s!(n-s)!}
\left(\pp_x^{n-s}\pp_p^sf\right)\left(\pp_x^s\pp_p^{n-s}g\right)=$$
$$=\sum_0^{\infty}\frac{\gk^n}{n!}\sum_0^n(-1)^r\left(\pp_x^r\pp_p^{n-r}f\right)
\left(\pp_x^{n-r}\pp_p^rg\right)$$
(note there are typos on p. 169 in \cite{ch}) and e.g. one has
\bq\label{2}
g*f=g(x+\gk\pp_p,p-\gk\pp_x)f=f(x-\gk\pp_p,p+\gd\pp_x)g
\end{equation}
The Moyal bracket can then be defined via
\bq\label{3}
\{f,g\}_M=\frac{1}{\gk}\{f\,Sin[\gk(\ola{\pp}_x\ora{\pp}_p-\ola{\pp}_p\ora{\pp}_x)]g\}
=\frac{1}{2\gk}(f*g-g*f)=
\end{equation}
$$=\sum_0^{\infty}\frac{(-1)^s\gk^{2s}}{(2s+1)!}\sum_0^{2s+1}(-1)^j\left(
\begin{array}{c}
2s+1\\
j\end{array}\right)[\pp_x^j\pp_p^{2s+1-j}f][\pp_x^{2s+1-j}\pp_p^jg]$$
corresponding to $\gk\to i\gk$ in \eqref{X18}.
\\[3mm]\indent
We emphasize also that many formulas in classical integrable systems already have a 
quantum mechanical (QM) flavor.  for example in \cite{ch,gal,sh} one shows how there
is a Moyal deformation $(KP)_M$ of dKP which for a particular value of $\gk$
($\gk =1/2$ in \cite{ch,gal}) creates an equivalence $(KP)_M\equiv (KP)_{Sato}$.
Actually QM features in integrable sysems seem inevitable because of Lax operator
formulations and the combinatorics inherent in Hirota equations and tau functions;
also early work by the Kyoto school provided many connections between KP and 
quantum field theory (QFT) (cf. \cite{cw}).  Such connections have since proliferated
in topological field theory (TFT), Seiberg-Witten (SW) theory, etc. where e.g.
effective actions can correspond to tau functions of integrable sysems and, somewhat
paradoxically, effective slow dynamics or Whitham dynamics (obtained by averaging out
fast fluctuations of angle variables) seems to correspond to a quantization (cf.
\cite{ch}, Chapter 5 or \cite{cag} for discussion).  On the other hand the so called
quantum inverse scattering method involving spin chains etc. for quantum integrable
systems (cf. \cite{ch,fu}), has a definite quantum group nature where the R-matrix
provides quasitriangularity.  The connection between R and r matrices leads one back to
classical dynamics but the theories for two types of integrable systems (classical and
quantum) have developed along different paths.  It seems that various discretizations
involving classical integrable systems (surveyed in \cite{ch}) should have a q-analysis
foundation and thus there may be other forms of
connecting glue between classical and quantum integrable systems
via discretization.  Indeed 
one almost seems to expect a discrete formulation to automatically have
quantum features.

\section{DISCRETIZATION AND MOYAL}
\renewcommand{\theequation}{2.\arabic{equation}}
\setcounter{equation}{0}

In \cite{ch} we expounded as some length on a series of papers by Kemmoku,
S. Saito, and collaborators (cf. \cite{ch} for references)
and we now want to organize some of this in a better manner and develop matters
somewhat further.  Thus we sketch first some fundamental ideas.  One defines 
\bq\label{6}
\nabla=\frac{e^{\gl\pp}-e^{-\gl\pp}}{2\gl}=\frac{1}{\gl}sinh(\gl\pp);\,\,\na_{{\bf a}}=
\frac{1}{\gl}sinh\left(\gl\sum a_i\pp_i\right)
\end{equation}
where $a_i\sim\pp/\pp_i$ and $\pp\sim\pp_x$.  Evidently ${\bf (A1)}\,\,\na
f(x)=(1/2\gl)[f(x+\gl)-f(x-\gl)]$ and $\na_{{\bf a}}f({\bf x})=(1/2\gl)[f({\bf x}
+{\bf a})-f({\bf x}-{\bf a})]$ (note that $\na_{{\bf a}}$ is not a vector).  Set then
(the
$a_i$ correspond to unspecified local coordinates $x_i$ generating a lattice with
vectors ${\bf a}$ in say ${\bf R}^N$ where $N\to \infty$ would require some
convergence stipulations)
\bq\label{7}
X^D=\int d{\bf a}v_{\gl}({\bf x},{\bf a})\na_{{\bf a}};\,\,\int d{\bf a}\sim\int\prod
da_i
\end{equation}
Next a difference one form is defined via ${\bf (A2)}\,\,\gO_D=\int d{\bf
a}w_{\gl}({\bf x},{\bf a})\gD^{{\bf a}}$ where $<\gD^{{\bf b}},\na_{{\bf a}}>=\gd({\bf
b}-{\bf a})$ and ($\vec{a}\sim{\bf a}$)
\bq
<\Omega^D,X^D>=\int d\vec{a}\int d\vec{b}<w_{\lambda}(\vec{x},\vec{b})
\Delta^{\vec{b}},v_{\lambda}(\vec{x},\vec{a})\nabla_{\vec{a}}>=
\int d\vec{a}w_{\lambda}(\vec{x},\vec{a})v_{\lambda}(\vec{x},\vec{a})
\label{XX20}
\end{equation}
Note also $\gD^{{\bf a}}$ can be realized via ($<\gD^{{\bf a}},\na_{{\bf b}}>=\gd({\bf
a}-{\bf b})$)
\bq
\gD^{{\bf a}}=\lambda csch[\lambda(\vec{a}\cdot\vec{\partial})]=\frac{2\lambda}
{e^{\lambda\vec{a}\cdot\vec{\partial}}-e^{-\lambda\vec{a}\cdot\vec
{\partial}}}=2\lambda\sum_0^{\infty}e^{-\lambda(2n+1)\vec{a}\cdot
\vec{\partial}}
\label{X21}
\end{equation}
In this connection we recall the $q^2$ difference operator ${\bf (A3)}\,\,
\pp_{q^2}f(x)=[f(q^2x)-f(x)]/[(q^2-1)x]$
with  ``dual" a Jackson integral ${\bf (A4)}\,\,\int_0^yd_{q^2}xf(x)=y(1-q^2)\sum_0^
{\infty}f(yq^{2n})q^{2n}$.
According to \cite{kas} there should be an unspecified
q-analysis version of \eqref{X21} related to pseudodifferential operators.  We can 
develop an interesting q-analysis counterpart to \eqref{X21} as follows.  Note first
that for $y=x+\gl$ one can write ${\bf
(A5)}\,\,(1/2\gl)[f(x+\gl)-f(x-\gl)]=[f(y+2\gl)-f(y)]/2\gl$ and for $q^2y=y+2\gl$ one
has $2\gl=(q^2-1)y$.  Then consider ${\bf (A6)}\,\,\tl{\na}=[exp(2\gl\pp)-1]/2\gl$
with 
\bq\label{10}
\tl{\na}f(y)=\frac{f(y+2\gl)-f(y)}{2\gl}=\frac{f(q^2y)-f(y)}{(q^2-1)y}=
\tl{\pp}_{q^2}f(y)\equiv
\end{equation}
$$\equiv \tl{\pp}_qf(z)=\frac{f(qz)-f(q^{-1}z)}{9q-q^{-1}z}\,\,\,(qy=z)$$
where $\tl{\pp}_{q^2}$ involves now a variable $q=q(y)$ if $\gl$ is to be regarded as
constant (alternatively one could regard $\gl$ as variable in $y$ and $q$ as  
constant or dispense with $\gl$ altogether).  For $\gl$ constant \eqref{X21} would
become formally a $y$ dependent inverse (note $(q^2-1)ny=2n\gl$)
\bq\label{11} 
\tl{\na}^{-1}=-2\gl(1-e^{2\gl\pp})^{-1}=-2\gl\sum_0^{\infty}
e^{2n\gl\pp}=(1-q^2)y\sum_0^{\infty} e^{(q^2-1)ny\pp}
\end{equation}
leading to 
\bq\label{12}
-2\gl(1-e^{2\gl\pp})^{-1}g(y)=G(y)=-2\gl\sum_0^{\infty} g(y+2n\gl)=
\end{equation}
$$=(1-q^2)y\sum_0^{\infty} g(y+(q^2-1)ny)$$
Evidently ${\bf (A7)}\,\,\tl{\na}G(y)=g(y)$  
so we can state (note a constant of
integration in \eqref{12} would vanish for $\int_0^yg\sim G(y)$)
\begin{proposition}
If we regard $q$ as $y$ dependent via $2\gl=(q^2-1)y$ with $\gl$ constant then the
inversion \eqref{12} has a modified Jackson type integral form
\bq\label{13}
G(y)={\ncint}_0^y\, g(x)d_{q^2}x\sim
-2\gl(1-e^{2\gl\pp})^{-1}g(y)=(1-q^2)y\sum_0^{\infty} g(y+(q^2-1)ny)
\end{equation}
\end{proposition}
{\bf REMARK 2.1.}
Note $y$ is fixed throughout so the calculations make sense and this reveals also a
property of Jackson integrals {\bf (A4)}, namely they do not seem to use the
integration variable $x$ at all (although change of variable techniques should work).
We emphasize that care is needed in using \eqref{10} in the form $\tl{\pp}_{q^2}$ when
computing $\tl{\pp}_{q^2}G(y)=g(y)$.  Thus $\tl{\pp}_{q^2}$ defined via $\tl{\na}$ in
\eqref{10} is \ul{{\bf not}} the same as $\pp_{q^2}$ unless provision is made for
$\gl=c$.  If we try to compute $\pp_{q^2}G(y)$ without keeping $\gl$ constant there
arises an awkward term $(1-q^2)q^2y\sum_0^{\infty}g(q^2y+(q^2-1)nq^2y)$ and 
$\pp_{q^2}G(y)\ne g(y)$.  The point is that $2\gl$ is constant and $(1-q^2)y=2\gl
\not\to (1-q^2)q^2y$.  Nor does $y+2n\gl=y+(1-q^2)ny$ go to $y+(1-q^2)nq^2y
=y+2n\gl q^2$
(rather e.g. $y+2n\gl\to q^2y+2n\gl=y+2(n+1)\gl=y+(1-q^2)(n+1)y$).  Thus for
$\tl{\pp}_{q^2} G(y)$ one must write
$(1/2\gl)[G(y+2\gl)-G(y)]=[(1-q^2)y]^{-1}[G(q^2y)-G(y)=\pp_{q^2} G(y)$ as desired.
If we regard this as a generally viable procedure of transferring ``standard"
differencing techniques in $\gl$ to q-analysis then constant $\gl$ steps for any y
correspond to constant steps $(1-q^2)y$ which means for large $y,\,\,q\to 1,$ so
if $G'$ is continuous for example then
\bq\label{14}
\tl{\pp}_{q^2}G(y)=\frac{G(q^2y)-G(y)}{(q^2-1)y}\sim\frac{g(y+2\gl)-
G(y)}{2\gl}=G'(\xi)
\end{equation}
for $y\leq \xi\leq y+2\gl=q^2y$ and for t large $y+2\gl\simeq y$ corresponds to
$q^2\to 1$.  There seems to be no reason not to use the $q,\,\gl$ correspondence in
general as long as computational consistency is maintained.
\\[3mm]\indent
{\bf REMARK 2.2.}  We will eventually dispense with $\gl$ altogether in rephrasing
matters entirely in $q$ so that $\tl{\pp}_{q^2}$ or $\tl{\pp}_q$ will not arise.
\\[3mm]\indent
Continuing now from \cite{ch} one can define difference 2-forms $\gO^D_2$, an exterior
difference operator $\gD$, and a Lie difference operator via (standard $\wedge$
product)
\bq\label{15}
\gO^D_2=\int d{\bf a}\int d{\bf b}w_{\gl}({\bf x},{\bf a},{\bf b})\gD^{{\bf
a}}\wg\gD^{{\bf b}};
\end{equation}
$$\gD\gO^D_2=\int d{\bf a}\int d{\bf b}\int d{\bf c}\na_{{\bf a}}w_{\gl}({\bf x},{\bf
a},{\bf b})\gD^{{\bf c}}\wg\gD^{{\bf a}}\wg\gD^{{\bf b}}$$
Since $[\na_{{\bf a}},\na_{{\bf b}}]=0$ one has $\gD\gD=0$ and finally for $X^D$ as in
\eqref{7}
\bq\label{16}
i_{\na_{{\bf c}}}(\gD^{{\bf a}}\wg\gD^{{\bf b}})=\gd({\bf c}-{\bf a})\gD^{{\bf b}}-
\gd({\bf c}-{\bf b})\gD^{{\bf a}};\,\,\mf{L}_{X^D}=\gD\cdot i_{X^D}+i_{X^D}\cdot\gD
\end{equation}
\indent
Now consider a phase space $\vec{x}\sim {\bf x}=(x,p)$ and in place of ${\bf
(A8)}\,\,X_fg=(f_p\pp_x-f_x\pp_p)g$ one writes ${\bf (A9)}\,\,X_f^D=\int
da_1da_2v_{\gl}[f](x,p,a_1,a_2)\na_{{\bf a}}$ where (cf. \eqref{7})
\bq
v_{\lambda}[f](x,p,a_1,a_2)=\left(\frac{\lambda}{2\pi}\right)^2\int db_1
db_2exp[-i\lambda(a_1b_2-a_2b_1)]f(x+\lambda b_1,p+\lambda b_2)
\label{XX2}
\end{equation}
which should correspond to $<\gD^{{\bf a}},X_f^D>$ (cf. Section 3).
Note $a_1b_2-a_2b_1$ can be written as $\vec{a}\times \vec{b}$ and
$(1/\gl)(\vec{a}\times\vec{b})$ is the area in $\gl$ units of the parallelogram formed
by $\vec{a}\times\vec{b}$ ($\gl$ is essentially a scaling factor here and not a 
Fourier variable).  
The symplectic structure of {\bf (A8)} is retained via an
interchange of $\vec{a}$ and $\vec{b}$.  We note that {\bf (A9)} can be written in the
form (the details are in \cite{ch})
\bq
X_f^D=\frac{-i\lambda}{(2\pi)^2}\int da_1da_2\int db_1db_2Sin
[\lambda(a_1b_2-a_2b_1)]f(x+\lambda b_1,p+\lambda b_2)e^
{\lambda(a_1\partial_x+a_2\partial_p)}
\label{XX24}
\end{equation}
leading to
\bq
X_f^Dg=
-\frac{-i\lambda}{(2\pi)^2}\int\int da_1da_2\int\int db_1db_2\times
\label{XX1}
\end{equation}
$$\times Sin[\lambda(a_1b_2-a_2b_1)]f(x+\lambda b_1,p+\lambda b_2)exp
[\lambda(a_1\partial_x+a_2\partial_p)]g(x,p)=$$
$$=-\frac{i\lambda}{(2\pi)^2}\int da\int db Sin[\lambda(a_1b_2-a_2b_1)]
f(x+\lambda b_1,p+\lambda b_2)g(x+\lambda a_1,p+\lambda a_2)$$
Subsequent calculation gives, using $x+\gl a_1=\ga_1$ and $p+\gl a_2=\ga_2$  (cf.
\cite{ch})
\bq
\int da\int db e^{i\lambda(a_1b_2-a_2b_1)}f(x+\lambda b_1,p+\lambda b_2)
g(x+\lambda a_1,p+\lambda a_2)=
\label{X230}
\end{equation}
$$=\frac{1}{\lambda^2}\left(\int f(x+i\lambda\partial_{\alpha_2},p-i\lambda
\partial_{\alpha_1})\int e^{i[b_2(\alpha_1-x)-b_1(\alpha_2-p)]}db\right)
g(\alpha_1,\alpha_2) d\alpha=$$
$$=\left(\frac{2\pi}{\lambda}\right)^2f(x-i\lambda
\partial_p,p+i\lambda\partial_x)g(x,p)\sim \left(\frac{2\pi}{\lambda}
\right)^2g*f$$
leading finally to
\bq\label{X24}
X_f^Dg=\frac{i}{\gl}Sin[\gl(\pp_{x_1}\pp_{p_2}-\pp_{p_1}\pp_{x_2})]
f(p_1,x_1)g(p_2,x_2)|_{(p,x)}=\{f,g\}_M
\end{equation}
In addition, from the Jacobi identity for the Moyal bracket one has
$$ [X^D_f,X^D_g]h=X^D_f\{g,h\}-X^D_g\{f,h\}=\{f,\{g,h\}\}-\{g,\{f,h\}\}\}=$$
\bq
=\{\{f,g\},h\}=X^D_{\{f,g\}}h
\label{X240}
\end{equation}
A symplectic form can also be given via
\bq\label{17}
\gO=\frac{1}{2\gl}\int\int da_1da_2\int\int db_1db_2e^{i\gl(a_1b_2-a_2b_1)}\gD^{{\bf
a}}\wg\gD^{{\bf b}}
\end{equation}
and this satisfies $i_{X_f^D}\gO=\gD f$ (analogous to $i_{X_f}\go=d\go$ for a
symplectic form $\go$).
Our formulas differ at times by $\pm i$ from \cite{kas,kaw} but everything seems
consistent and correct here;
the philosophy of running $a_i$ over ${\bf R}\sim (-\infty,\infty)$ is crucial in the
calculations (alternatively $\int$ could represent a sum over a discrete symmetric
set, e.g. $[-N,N]$ with N infinite or not).  We note also a somewhat quasi Fourier
theoretic version of the formulas {\bf (A9)}, \eqref{XX2}, \eqref{XX24}, etc.  
developed in \cite{ch}.  Thus
consider
\bq
v_{\lambda}[f]({\bf x},{\bf a})=\left(\frac{\lambda}{2\pi}\right)^2\int
d{\bf b}e^{-i\lambda({\bf a}\times {\bf b})}e^{\lambda \vec{b}\cdot
\vec{\partial}}f
\label{120J}
\end{equation}
Hence (using ${\bf b}\to-{\bf b}$)
\bq
v_{\lambda}[f]({\bf x},-{\bf a})=\left(\frac{\lambda}{2\pi}\right)^2
\int d{\bf b}e^{i\lambda({\bf a}\times {\bf b})}e^{\lambda\vec{b}\cdot\vec{\partial}}f=
\left(\frac{\lambda}{2\pi}\right)^2\int d{\bf b}e^{-i\lambda({\bf a}\times
{\bf b})}e^{-\lambda\vec{b}\cdot\vec{\partial}}f
\label{121J}
\end{equation}
and since $\nabla_{-{\bf a}}=-\nabla_{{\bf a}}$ one gets
\bq
X^D_f=\int d{\bf a}v_{\lambda}[f]({\bf x},{\bf a})\nabla_{{\bf a}}=-\int_{\infty}^
{-\infty}d{\bf a}v_{\lambda}[f]({\bf x},-{\bf a})\nabla_{-{\bf a}}=-\int
d{\bf a}v_{\lambda}[f]({\bf x},-{\bf a})\nabla_{{\bf a}}
\label{122J}
\end{equation}
Consequently
\bq
X^D_f=\frac{1}{2}\int d{\bf a}\left[v_{\lambda}
[f]({\bf x},{\bf a})-v_{\lambda}[f]
({\bf x},-{\bf a})\right]\nabla_{{\bf a}}=
\label{123J}
\end{equation}
$$=\frac{\lambda^3}{4\pi^2}\int d{\bf a}\int d{\bf b}e^{-i\lambda({\bf a}
\times {\bf b})}\left\{\frac{e^{\lambda\vec{b}\cdot\vec{\partial}}-e^{-
\lambda\vec{b}\cdot\vec{\partial}}}{2\lambda}\right\}f\nabla_{{\bf a}}=
\frac{\lambda^3}{4\pi^2}\int d{\bf a}\int d{\bf b}e^{-i\lambda({\bf a}
\times {\bf b})}\nabla_{{\bf b}}f\nabla_{{\bf a}}$$
This formula provides another representation for $X^D_f$ via
\bq
X^D_f=\int d{\bf a}\tilde{v}_{\lambda}[f]({\bf x},{\bf a})\nabla_{{\bf a}};\,\,
\tilde{v}_{\lambda}[f]({\bf x},{\bf a})=\frac{\lambda^3}{4\pi^2}\int d{\bf b}
e^{-i\lambda({\bf a}\times {\bf b})}\nabla_{{\bf b}}f
\label{125J}
\end{equation}
\indent
The above gives a direct discretization of phase space and the natural difference
analogue of Lie bracket leads to the Moyal bracket.  Thus one takes $\gl\sim \hbar/2$
and defines $X_A^Q=\hbar X_A^D$ for functions $A(x,p)$ and there is a Heisenberg
equation ($H\sim$ Hamiltonian) ${\bf (A10)}\,\,-i\hbar\pp_tX_A^Q=[X_A^Q,X_H^Q]$ (where
both $A$ and $H$ may contain $\hbar$).  This is compatible with ${\bf (A11)}\,\,\pp_tA=
\{A,H\}_M$ (cf. \eqref{X24}, \eqref{X240}).  To see how this works we recall the
standard quantum mechanical (QM) idea of Wigner distribution function $F_w$ with
$\int F_wdx=1$ and $<\hat{A}>=\int F_wAdx$ for he expectation value of an operator
$\hat{A}$ associated to the observable function A (Weyl ordering is to be invoked when
ordering is needed and details are in \cite{ch}).
The corresponding discrete version is given via a difference 1-form
\bq\label{X28}
P_{F_w}=\frac{\hbar}{4}\int\int da_1da_2\int\int db_1db_2e^{i\hbar(a_1b_2-
a_2b_1)/2}F_w\left(x+\frac{\hbar}{2}b_1,p+\frac{\hbar}{2}b_2\right)\Delta^a
\end{equation}
so ${\bf (A12)}\,\,<P_{F_w},X_A^Q>=\int dxdp\,F_w(x,p)A(x,p)=<\hat{A}>$.  In the
Heisenberg picture the time dependence is ${\bf
(A13)}\,\,\pp_t<P_{F_w},X_A^Q>=<P_{F_w},X_A^Q(t)>$ which in the Schr\"odinger picture
becomes ${\bf (A14)}\,\,\pp_t<P_{F_w},X_A^Q>=<P_{F_w}(t),X_A^Q>$.  Here the solution
of  {\bf (A10)} necessarily is 
\bq\label{18}
X_A^Q(t)=exp\left(-\frac{it}{\hbar}X_H^Q\right)X_A^Qexp\left(\frac{it}{\hbar}
X_H^Q\right)
\end{equation}
(simply differentiate $X_A^Q=exp[(it/\hbar)X_H^Q]X_A^Q(t)exp[(-it/\hbar)X_H^Q]$ and
note that in {\bf (A10)} $X_A^Q\sim X_A^Q(t)$).  This corresponds to a solution of 
{\bf (A11)} of the form ${\bf (A15)}\,\,A(t)=[exp(it/\hbar)X_H^Q]A$ and in the
Heisenberg picture
\bq\label{19}
-i\hbar\frac{d}{dt}<P_{F_w},X^Q_A(t)>=<P_{F_w},[X^Q_A(t),X^Q_H]>=
<P_{F_w},X^Q_{\{A(t),H\}_M}>
\end{equation}
where the right side is $<P_{\{H,F_w(t)\}_M},X_A^Q>$ upon defining ${\bf (A16)}
\,\,F_w(t)=exp[-(it/\hbar)X_H^Q]F_w$ so that ${\bf (A17)}\,\,\pp_tP_{F_w(t)}=
P_{\{H,F_w(t)\}_M}\equiv \pp_tF_w(t)=\{H,F_w(t)\}_M$.

\section{Q-DISCRETIZATION}
\renewcommand{\theequation}{3.\arabic{equation}}
\setcounter{equation}{0}

Let us consider now a variation on Section 2 based on a q-lattice.  This will
constitute a different approach from those in Remark 2.1 and Proposition 2.1 in that
we keep q fixed.  Indeed q can play the role of $\gl$ and we write
\bq\label{20}
\hat{\na}_{mn}f(x,p)=\frac{f(xq^{2m},pq^{2n})-f(x,p)}{(q^{2m}-1)x(q^{2n}-1)p}
\end{equation}
\bq\label{30}
\check{\na}_{mn}g(x,p)=\frac{g(xq^m,pq^n)-g(xq^{-m},pq^{-n})}{(q^m-q^{-m})(q^n-q^{-n})xp}=
\end{equation}
$$\frac{e^{\gl(m,n)\cdot(\hat{\pp}_1,\hat{\pp}_2)}-e^{-\gl(m,n)
\cdot(\hat{\pp}_1,\hat{\pp}_2)}}{(q^m-q^{-m})(q^n-q^{-n})xp}\hat{g}(log(x),log(p))
=q^{m+n}e^{-\gl(m,n)\cdot(\hat{\pp}_1,\hat{\pp}_2)}\hat{\na}_{mn}g=G$$
so $(m,n)$ plays the role of Fourier variables$(a_1,a_2)\sim{\bf a}$.  We recall from
\cite{ch} the device ${\bf
(A18)}\,\,\gl=log(q),\,\,exp(\gl)=q,\,\,f(x)=\hat{f}(log(x)),\,\,q^{2mx\pp_x}f(x)=
exp[2m\gl\pp_{log(x)}]\hat{f}(log(x))=\hat{f}(log(x)+2mlog(q))=\hat{f}(log(q^{2m}x))=
f(xq^{2m})$.  This suggests an inversion for $\hat{\na}_{mn}$ written via
\bq\label{21}
\hat{\na}_{mn}f(x,p)=\frac{(e^{2\gl (m,n)\cdot(\hat{\pp}_1,\hat{\pp}_2)}-1)}
{(q^{2m}-1)x(q^{2n}-1)p}\hat{f}(log(x),log(p))
\end{equation}
($\hat{\pp}_1=\pp/\pp\,log(x),\,\,\hat{\pp}_2=\pp/\pp\,log(p)$) in a form similar to a 
Jackson integral.  Thus first we can derive a Jackson integral as follows.  Write
\bq\label{22}
\na f(x)=\pp_{q^2}f(x)=\frac{f(q^2x))-f(x)}{(q^2-1)x}=\frac{(e^{2\gl
x\pp_x}-1)}{(q^2-1)x}f(x)=g(x)
\end{equation}
with formally
\bq\label{23} 
f(x)=(1-q^2)\sum_0^{\infty}e^{2k\gl x\pp_x}(xg(x))=(1-q^2)\sum_0^{\infty}q^{2k}x
g(q^{2k}x)
\end{equation}
which is the Jackson integral $\int_0^xd_{q^2}yg(y)$.  Similarly we can write now
formally
\bq\label{24}
\hat{\na}_{mn}^{-1}g(x,p)=-(q^{2m}-1)(q^{2n}-1)\sum_0^{\infty}e^{2\gl k(m,n)\cdot(\hat
{\pp}_1,\hat{\pp}_2)}(xpg(x,p))=
\end{equation}
$$=-(q^{2m}-1)(q^{2n}-1)\sum_0^{\infty}q^{2mk}xq^{2nk}pg(q^{2mk}x,q^{2nk}p)=G(x,p)$$
This can be checked via
\bq\label{25}
\frac{G(q^{2m}x,q^{2n}p)-G(x,p)}{(q^{2m}-1)x(q^{2n}-1)p}=g(x,p)=
\end{equation}
$$=-\sum_0^{\infty}q^{2m(k+1)}q^{2n(k+1)}g(q^{2m(k+1)}x,q^{2n(k+1)}p)+\sum_0^{\infty}
q^{2mk}q^{2nk}g(q^{2mk}x,q^{2nk}p)$$
Hence we have proved
\begin{proposition}
The difference operator $\hat{\na}_{mn}$ of \eqref{20} can be inverted via \eqref{24}
s a kind of extended Jackson integral.  Similarly one has
\bq\label{31}
\check{\na}_{mn}^{-1}g(x,p)=q^{-m-n}\hat{\na}^{-1}_{mn}g(xq^{-m},pq^{-n})=
\end{equation}
$$=-q^{-m-n}(q^{2m}
-1)(q^{2n}-1)\sum_0^{\infty}q^{2mk-m}xq^{2nk-n}pg(q^{2mk-m}x,q^{2nk-n}p)=$$
$$=-(q^m-q^{-m})(q^n-q^{-n})xp\sum_0^{\infty}q^{(2k-1)(m+n)}g(q^{(2k-1)m}x,
q^{(2k-1)n}p)$$
\end{proposition}
It should be possible now to duplicate most of the machinery in Section 2 with q
discretization as above.  We note that this 
procedure and the resulting formulas appear to be different from any of the
phase space discretizations in \cite{cze,dae,fe,he,kzq,lh,sz,sza,smi,tb,wzz}.
We will consider an analogue of $X_f^D$ in {\bf (A9)} or \eqref{125J} via
\bq\label{26}
\hat{X}_f^D=\sum_{m,n}v_q[f](x,p,m,n)\hat{\na}_{mn}\,\,{\bf or}\,\,\check{X}_f^D=\sum
v_q[f](x,p,m,n)\check{\na}_{mn}
\end{equation}
where we need then a formula for $v_q[f]$ which can perhaps be modeled on \eqref{125J}
in a quasi Fourier spirit.  Note that the stipulation $<\gD^{{\bf a}},\na_{{\bf
b}}>=\gd({\bf a}-{\bf b})$, or $\hat{\gD}^{mn}=\hat{\na}^{-1}_{mn}$ as in 
\eqref{24}-\eqref{25} simply provides a tautology ${\bf
(A19)}\,\,v_q[f](x,p,m,n)=<\hat{\gD}^{mn},\hat{X}_f^D>$ or as in \eqref{XX2} the
equation ${\bf (A20)}\,\,<\gD^{{\bf a}},X_f^D>=<\gD^{{\bf a}},\int d{\bf
b}v_{\gl}[f](x,p,{\bf b})\na_{{\bf b}}>=v_{{\bf a}}[f](x,p,{\bf a})$.
Thus one should realize that $v_{\gl}[f]$ is simply selected in an ad hoc manner so
that $X_f^Dg=\{f,g\}_M$.  
It turns out that the use of $\hat{\na}_{mn}$ and $\hat{X}_f^D$ would not reproduce
a suitable $\pm$ symmetry for a quasi Fourier approach so we will concentrate on
$\check{X}_f^D$ and $\check{\na}_{mn}$.
\\[3mm]\indent
In \cite{dam} a quantum q-Moyal bracket 
$(\hbar\ne 0$) is suggested in the form
\bq\label{288}
\{p^mx^n,p^kx^{\ell}\}_{qM}=\frac{1}{i\hbar}(q^{nk}p^mx^n *
p^kx^{\ell}-q^{m\ell}p^kx^{\ell}*p^mx^n)
\end{equation}
where $*$ can refer to standard or antistandard orderings via
($\nu=log(q)$ and $D_z\sim\pp_q$)
\bq\label{289}
*_S\equiv\sum_0^{\infty}\frac{(i\hbar)^r}{[r]!}\ola{D}^r_pexp(\nu\ola{\pp}_ppx
\ora{\pp}_x)\ora{D}^r_x;
\end{equation}
$$*_A\equiv\sum_{s=0}^{\infty}(-\nu\ola{\pp}_xx)^s\sum_{r=0}^s
\frac{(-i\hbar)^rq^{r(r-1)/2}}{[r]!}\ola{D}_x^r\ora{D}_p^r(p\ora{\pp}_p)^s$$
Here standard ordering involves XP products and antistandard has PX products
(see Section 5).  The symbol map is $S_S(X^mP^n)=S_A(P^mX^n)=p^mx^n$; Weyl ordering
is also considered but there are some complications.  We note also for $\hbar=0$ one
has classical star products based on ($\nu=log(q)$ - cf. \cite{dam,dzz})
\bq\label{290}
*^q_S\equiv exp(\nu\ola{\pp}_ppx\ora{\pp}_x);\,\,*_A^q\equiv
exp(-\nu\ola{\pp}_xxp\ora{\pp}_p);
\end{equation}
$$*_W^q\equiv
exp\left(-\frac{\nu}{2}(\ola{\pp}_xxp\ora{\pp}_p)-\ola{\pp}_ppx\ora{\pp}_x)\right)$$
(here $*_W^q$ refers to Weyl ordering); these star products all satisfy
\bq\label{291}
q^{nk}p^mx^n*^qp^kx^{\ell}-q^{m\ell}p^kx^{\ell}*^qp^mx^n=0
\end{equation}

\section{CALCULATIONS}
\renewcommand{\theequation}{4.\arabic{equation}}
\setcounter{equation}{0}

For completeness we will give a number of calculations to show how our results are
parallel to Section 2 and can be reached through some quasi Fourier type procedures.
First we recall some useful formulas (cf. \cite{cw,fz,kca}), namely
\bq\label{35}
\gd(z-w)=z^{-1}\sum_{n\in{\bf Z}}\left(\frac{z}{w}\right)^n=z^{-1}\hat{\gd}(q/w)
\end{equation}
There are many nice calculations available using \eqref{35}; we mention e.g.
($Res_z\sum a_nz^n=a_{-1}$ and $D_z=z(d/dz)$)
\bq\label{36}
\gd(w-z)=w^{-1}\sum_{{\bf Z}}\left(\frac{w}{z}\right)=w^{-1}\sum_{{\bf
Z}}\left(\frac{z}{w}\right)^n= z^{-1}\sum\left(\frac{z}{w}\right)^n=\gd(z-w);
\end{equation}
$$Res_zf(z)\gd(z-w)=f(w);\,\,f(z)\hat{\gd}(az)=f(a^{-1})\hat{\gd}(az);\,\,Res_z\pp
a(z)b(z)=
-Res_za(z)\pp b(z)$$
This will provide a delta function corresponding to
$\int exp[ib_2(\ga_1-x)-ib_1(\ga_2-p)]d{\bf b}$.  Now, leaving aside possible
multiplicative factors (cf. Remark 4.1),
consider \eqref{XX2} in the form 
\bq\label{29}
v_q[f](x,p,{\bf a})=c(q)\sum_{r,s}q^{ms-nr}f(q^rx,q^sp)
\end{equation} 
leading to
(cf. \eqref{XX24} - \eqref{XX1})
\bq\label{32} 
X_f^D=\hat{c}(q)\sum_{m,n,r,s}(q^{ms-nr}-q^{-ms+nr})
f(q^rx,q^sp)\cdot
q^{(m,n)\cdot(\hat{\pp}_1,\hat{\pp}_2)};
\end{equation}
$$X_f^Dg=\hat{c}(q)\sum_{m,n,r,s}(q^{ms-nr}-q^{-ms+nr})
f(q^rx,q^sp,)g(q^mx,q^np)$$
while \eqref{X230} can be written as ($x+\lambda a_1=\alpha_1$
and $p+\lambda a_2=\alpha_2$)
\bq
\int da\int db e^{i\lambda(a_1b_2-a_2b_1)}f(x+\lambda b_1,p+\lambda b_2)
g(x+\lambda a_1,p+\lambda a_2)
\label{33}
\end{equation}
$$\frac{1}{\lambda^2}\int\int d\alpha db e^{i[b_2(\alpha_1-x)-b_1(\alpha_2-p)]}
f(x+\lambda b_1,p+\lambda b_2)g(\alpha_1,\alpha_2)=$$
$$=\frac{1}{\lambda^2}\left(\int f(x+i\lambda\partial_{\alpha_2},p-i\lambda
\partial_{\alpha_1})\int e^{i[b_2(\alpha_1-x)-b_1(\alpha_2-p)]}db\right)
g(\alpha_1,\alpha_2) d\alpha=$$
$$=\left(\frac{2\pi}{\lambda}\right)^2\int \left[f(x+i\lambda\partial_{\alpha_2},
p-i\lambda\partial_{\alpha_1})\delta(\alpha_1-x,\alpha_2-p) \right]g(\alpha_1,
\alpha_2)d\alpha=$$
$$=\left(\frac{2\pi}{\lambda}\right)^2f(x-i\lambda
\partial_p,p+i\lambda\partial_x)g(x,p)\sim \left(\frac{2\pi}{\lambda}
\right)^2g*f$$
Intuitively one thinks of $\gl\sim log(q)$, ${\bf a}\sim (m,n)$,
and ${\bf b}\sim (r,s)$ so the substitution
$x+\gl a_1=\ga_1$ corresponds to $\ga_1/x=q^m$; similarly $\ga_2/p=q^n$ and the second
and third lines in \eqref{33} correspond to
\bq\label{34}
\gG_1=c(q,p,x)\sum_{\ga}\sum_{r,s}\left(\frac{\ga_1}{x}\right)^s\left(\frac{\ga_2}{p}\right)^r
f(xq^r,pq^s)g(\ga_1,\ga_2)
\end{equation}
where $\sum_{\ga}\sim Res_{{\bf \ga}}(1/\ga_1\ga_2)$.  The first question is to ask
if we can write something like
\bq\label{39}
\sum_{r,s}f(xq^r,pq^s)\left(\frac{\ga_1}{x}\right)^s\left(\frac{\ga_2}{p}\right)^{-r}
\sim
f(xq^{\hat{\pp}_1},pq^{-\hat{\pp}_2})\hat{\gd}\left(\frac{\ga_1}{x}\right)\hat{\gd}
\left(\frac{p}{\ga_2}\right)
\end{equation}
in analogy to lines 3 and 4 of \eqref{33}.  We could imagine e.g. $f(x,p)=\sum
a_{k\ell}x^kp^{\ell}$ and look at
$$\sum_{r,s}x^kp^{\ell}q^{kr}q^{\ell
s}\left(\frac{\ga_1}{x}\right)^s\left(\frac{\ga_2}{p}\right)^{-r}=\sum_{r,s}
x^kp^{\ell}q^{-k\hat{\pp}_2}q^{\ell\hat{\pp}_1}\left(\frac{\ga_1}{x}\right)^s
\left(\frac{\ga_2}{p}\right)^{-r}=$$
\bq\label{40}
=x^kp^{\ell}q^{-k\hat{\pp}_2}q^{\ell\hat{\pp}_1}\hat{\gd}\left(\frac{\ga_1}{x}\right)
\hat{\gd}\left(\frac{p}{\ga_2}\right)
\end{equation}
since $q^{-k\hat{\pp}_2}(\ga_2/p)^{-r}=(q^{-k}\ga_2/p)^{-r}=q^{kr}(\ga_2/p)^r$.
Consequently for $f=\sum a_{k\ell}x^kp^{\ell}$ in \eqref{34} we have
\bq\label{41}
\gG_1=c(q,p,x)\sum_{k,\ell}a_{k\ell}x^kp^{\ell}q^{-kp\pp_p}q^{\ell
x\pp_x}g(x,p)
\end{equation}
since $Res_{{\bf \ga}}(1/\ga_1\ga_2)\hat{\gd}(\ga_1/x)\hat{\gd}(p/\ga_2)g(\ga_1,\ga_2)
=g(x,p)$ and e.g. $\hat{\pp}_1$ in $\ga_1$ becomes $\hat{\pp}_1=x\pp_x$.  This leads to
\bq\label{42}
\gG_1=c(q,p,x)\sum a_{k\ell}x^kp^{\ell}g(xq^{\ell},pq^{-k})
\end{equation}
as a putative $g*f$ (cf. \eqref{33}).  For $g=\sum b_{\gag\gb}x^{\gag}p^{\gb}$ 
this corresponds to
\bq\label{43}
\gG_1(f,g)=c\sum_{k,\ell,\gag\gb}a_{k\ell}b_{\gag\gb}x^{k+\gag}p^{\ell+\gb}
q^{\ell\gag-k\gb}\sim g*f
\end{equation}
The terms of the form \eqref{X230} corresponding to $exp[-i\gl(a_1b_2-a_2b_1)]$ in 
\eqref{XX1} involve now in place of \eqref{39} a term
\bq\label{44}
-\sum_{r,s}f(xq^r,pq^s)\left(\frac{\ga_1}{x}\right)^{-s}\left(\frac
{\ga_2}{p}\right)^r=-f(xq^{-\hat{\pp}_1},pq^{\hat{\pp}_2})\hat{\gd}\left(\frac{x}
{\ga_1}\right)\hat{\gd}\left(\frac{\ga_2}{p}\right)
\end{equation}
Hence we get for $f$ and $g$ as before
\bq\label{46}
\gG_2\sim f*g=-c(q,p,x)\sum
a_{k\ell}x^kp^{\ell}b_{\gag\gb}(xq^{-\ell})^{\gag}(pq^k)^{\gb}=c\sum a_{kl}
b_{\gag\gb}x^{k+\ga}p^{\ell+\gb}q^{k\gb-\ell\gag}
\end{equation}
leading to 
\begin{proposition}
For $f(x,p)=\sum a_{k\ell}x^kp^{\ell}$ and $g(x,p)=\sum b_{\gag,\gb}x^{\gag}p^{\gb}$
one obtains in an heuristic manner
\bq\label{47}
\{f,g\}_M\sim f(xq^{-p\pp_p},pq^{x\pp_x})g(x,p)-g(xq^{-p\pp_p},pq^{x\pp_x})f(x,p)\sim
\end{equation}
$$\sim c(q,p,x)\sum_{k,\ell,\gag,\gb}a_{k\ell}b_{\gag\gb}x^{k+\gag}
p^{\ell+\gb}\left(q^{k\gb-\ell\gag}-q^{\ell\gag-k\gb}\right)$$
where $c(q,p,x)$ is to be stipulated (cf. Corollary 4.1 for an essentially
equivalent formula).
Note by inspection or construction $\{f,g\}_M=-\{g,f\}_M$.
\end{proposition}
\indent
If we use the formulation of \eqref{120J} - \eqref{125J} a slightly different formula
emerges involving a multiplicative factor which is missed by the analogy constructions
above.
Thus we check the passage \eqref{120J} to \eqref{125J}.
\eqref{120J} is the same as \eqref{XX2} corresponding to \eqref{29} and \eqref{121J}
corresponds to
\bq\label{48}
v_q[f](x,p,-{\bf a})\sim c\sum_{r,s}q^{ms-nr}f(q^{-r}x,q^{-s}p)
\end{equation}
which would follow from \eqref{29} by sending $(m,n)\to -(m,n)$ and $(r,s)\to -(r,s)$.
This makes sense if the sums are $-\infty\to\infty$ and there seems to be no objection
to that.  Then one would have (taking now $\na_{{\bf a}}\sim\check{\na}_{mn}$ as in
\eqref{30})
\bq\label{49}
X_f^D=\int d{\bf a}v_{\gl}[f]({\bf x},{\bf a})\na_{{\bf a}}\sim
c\sum_{m,n}v_q[f](x,p,m,n)\check{\na}_{mn}=
\end{equation}
$$=c\sum_{m,n}\sum_{r,s}q^{ms-nr}f(q^rx,q^sp)\frac{q^{mx\pp_x}q^{np\pp_p}-q^{-mx\pp_x}
q^{-np\pp_p}}{(q^m-q^{-m})(q^n-q^{-n})xp}$$
\bq\label{50}
X_f^D=-\int d{\bf a}v_{\gl}[f]({\bf x},{\bf -a})\na_{{\bf a}}\sim X_f^D=c\sum_{m,n}
v_q[f](x,p,-m,-n)\check{\na}_{-m,-n}=
\end{equation}
$$=c\sum_{m,n,r,s}q^{ms-nr}f(q^{-r}x,q^{-s}p)\frac{q^{-mx\pp_x}q^{-np\pp_p}-q^{mx\pp_x}
q^{np\pp_p}}{(q^{-m}-q^m)(q^{-n}-q^n)xp}=$$
$$=-c\sum_{m,n,r,s}q^{ms-nr}f(q^{-r}x,q^{-s}p)\frac{q^{mx\pp_x}q^{np\pp_p}-
q^{-mx\pp_x}q^{-np\pp_p}}{(q^m-q^{-m})(q^n-q^{-n})xp}$$
exactly as in \eqref{122J} (note the minus sign appears in the last equation instead
of at the beginning).  Hence
$$X_f^D=\frac{1}{2}(\eqref{49}+\eqref{50})=c\sum_{m,n,r,s}q^{ms-nr}
[f(q^rx,q^sp)-f(q^{-r}x,q^{-s}p)]\check{\na}_{mn}=$$
\bq\label{51}
= c\sum_{m,n,r,s}q^{ms-nr}(q^r-q^{-r})(q^s-q^{-s})xp\check{\na}_{rs}f\check{\na}_{mn}
\end{equation}
which is a difference version of \eqref{123J}.  One sees that factors of
$(q^r-q^{-r}),\,\,(q^s-q^{-s}),\,\,(q^m-q^{-m})$, and $(q^n-q^{-n})$ have become
involved in place of powers of $\gl$ and this must be clarified; otherwise the
patterns go over.
\\[3mm]\indent
To clarify we compare \eqref{32} and \eqref{51} and write \eqref{51} in the form
\bq\label{52}
{}_1X_f^Dg=
\end{equation}
$$=\frac{c}{xp}\sum q^{ms-nr}[f(q^rx,q^sp)-f(q^{-r}x,q^{-s}p)]G(q,m,n)[g(xq^m,pq^n)-
g(xq^{-m},pq^{-n})]$$
where $G^{-1}(q,m,n)=(q^m-q^{-m})(q^n-q^{-n})=G(q,-m,-n)$.  Set $f_{\pm}\sim 
f(q^{\pm m}x,q^{\pm n}p)$ so in an obvious notation
\bq\label{53}
{}_1X_f^Dg=\frac{c}{xp}\sum
q^{ms-nr}G(q,m,n)[f_{+}g_{+}+f_{-}g_{-}-f_{+}g_{-}-f_{-}g_{+}]
\end{equation}
Now evidently, changing $m,n\to -m,-n$, one obtains a formula $\sum
q^{ms-nr}G(q,m,n)f_{+}g_{-}\to\sum e^{-ms+nr}G(q,m,n)f_{-}g_{+}$, etc. so
\bq\label{54}
{}_1X_f^Dg=\frac{c}{xp}\sum q^{-ms+nr}(f_{+}g_{-}+f_{-}g_{+}-f_{+}g_{+} - f_{-}g_{-}
\end{equation}
leading to
\bq\label{55}
{}_1X_f^Dg=cxp\sum
G(q,m,n)\times
\end{equation}
$$\times\left(q^{ms-nr}-q^{-ms+nr}\right)[f(q^rx,q^sp)-f(q^{-r}x,q^{-s}p)]
[g(xq^m,pq^n)-g(xq^{-m},pq^{-n})]$$
This is similar to \eqref{32} which has the form
\bq\label{56}{}_2X_f^Eg=\hat{c}\sum\left(q^{ms-nr}-q^{-ms+nr}\right)f_{+}g_{+}
\end{equation}
$$=\hat{c}\sum (\,\,\,)f_{-}g_{-}=-\hat{c}\sum (\,\,\,)f_{+}g_{-}=-\hat{c}\sum
(\,\,\,)f_{-}g_{+}$$
which implies
\bq\label{57}
{}_2X_f^Dg=
\end{equation}
$$=\hat{c}\sum\left(q^{ms-nr}-q^{-ms+nr}\right)[f(q^rx,q^sp)-f(q^{-r}x,q^{-s}p)]
[g(xq^m,pq^n)-g(xq^{-m},pq^{-n})]$$
This is essentially the same as ${}_1X_f^D$ except for the $G(q,m,n)$ factor.
For esthetic reasons one prefers the form ${}_1X_f^Dg$ since it has the more visibly
meaningful form \eqref{49} and $\gl$ plays a consistent role (cf. Remark 4.1 below).
Thus in summary
\begin{proposition}
The difference version of Section 2 can be expressed via
\bq\label{58}
X_f^D=\sum_{m,n}v_q[f](x,p,m,n)\check{\na}_{mn};\,\,v_q[f]=\sum_{r,s}
q^{ms-nr}f(q^rs,q^sp);
\end{equation}
$$\check{\na}_{mn}g=\frac{g(xq^m,pq^n)-g(xq^{-m},pq^{-n})}
{(q^m-q^{-m})(q^n-q^{n})xp};$$
$$X_f^Dg=\frac{1}{2xp}\sum_{m,n,r,s}q^{ms-nr}\frac{[f(q^rx,q^sp)-f(q^{-r}x,q^{-s}p)]
[g(q^mx,q^np)-g(q^{-m}x,q^{-n}p)]}{(q^m-q^{-m})(q^n-q^{-n})}$$
\end{proposition}
The latter expression is our putative Moyal bracket and one has
\begin{corollary}
Writing out $X_d^Dg$ for monomials $f=x^ap^b$ and $g=x^cp^d$ yields
\bq\label{59}
X_f^Dg=\{f,g\}_M=\frac{1}{2xp}\sum q^{ms-nr}\frac{x^{a+c}p^{b+d}[(q^{ra+bs}-q^{-ra-bs})
(q^{mc+nd}-q^{-mc-nd})]}{(q^m-q^{-m})(q^n-q^{-n})}
\end{equation}
Further since, as in \eqref{53} - \eqref{55}, one has
$-\sum_{m,n}q^{ms-nr}Gg_{-}=-\sum_{m,n}q^{-ms+nr}Gg_{+}$ and
$-\sum_{r,s}q^{ms-nr}f_{-}=-\sum q^{-ms+nr}f_{+}$ there results
\bq\label{60}
X_f^Dg=\frac{1}{2xp}\sum\frac{q^{ms-nr}-q^{-ms+nr}}{(q^m-q^{-m})(q^n-q^{-n})}
f(q^rx,q^sp)g(q^mx,q^np)
\end{equation}
This is reminiscent of \eqref{32} but with a $G(q,m,n)$ factor so the
calculation \eqref{47} applies with $G(q,m,n)$ inserted and consequently 
$\{f,g\}_M=-\{g,f\}_M$ as before, although this is not immediately visible from
\eqref{60}.  We note also from \eqref{58} or \eqref{60} that it does no harm to use
alternatively a form based on \eqref{125J} in the form (cf. \eqref{51})
\bq\label{61}
X_f^Dg=c\sum_{m,n,r,s} q^{ms-nr}\check{\na}_{rs}f\check{\na}_{mn}g
\end{equation}
which inserts an additional factor $G(q,r,s)$ into \eqref{58}.
\end{corollary}
\indent
{\bf REMARK 4.1.}
The multiplicative factors involve terms $(q^m-q^{-m}),\,\,(q^n-q^{-n}),\,\,
(q^r-q^{-r}),$ or $(q^s-q^{-s})$, all of which correspond to a $\gl$ arising from
$\check{\na}_{mn}$ or $\check{\na}_{rs}$; instead of coming out of the integral signs
as $\gl$ in the continuous versions of Section 2 they have to be summed.  Note the 
correspondence $x+\gl a_1=\ga_1$ corresponding to $\ga_1/x=q^m$ uses $\gl$ in a
different manner so it is at first glance surprising that ${}_2X_f^D$ even comes close
to ${}_1X_f^D$.  The relations of our formulas to the star products and Moyal brackets
of \eqref{288} - \eqref {291} will be examined later as well as the expansion of
material in \cite{ch} related to work of Curtright, Fairlie, Zachos, and the Saito
school (cf. \cite{ch} for references).  We note also that for a complex phase space
$\{z,\gz\}$ (not clarified) an 
interesting variation on the q-Moyal bracket of \eqref{47} or
\eqref{58} is given in \cite{kaw} for a KP situation
(cf. also \cite{ch} where this is expanded).  This is applied to a KP hierarchy
context using complex variable methods and, although powers of q are inserted
in various places, it is not developed systematically in a q-analysis manner and
no recourse to q-derivatives is indicated.  We will expand further the treatment of
\cite{ch} for this situation in a subsequent paper.

\end{document}